\documentclass[12pt,twoside,reqno]{amsproc}
\usepackage[centertags]{amsmath}
\usepackage{amssymb,amsfonts,amsthm}
\usepackage{pictex}

\newcommand{\TryPackage}[3]{\IfFileExists{#1.sty}{\usepackage{#1}
#2}{#3}}
\TryPackage{mathrsfs}{\renewcommand{\mathcal}{\mathscr
}}{%
    \TryPackage{eucal}{}{}}
\usepackage{a4mlnew,math401,local} 

\begin{document}

\title{Unbounded Fredholm Operators and Spectral Flow}

\author{Bernhelm Booss--Bavnbek}
\address{Institut for Matematik og Fysik\\ 
Roskilde Universitetscenter\\
DK--4000 Roskilde\\
Denmark}
\email{booss@mmf.ruc.dk}
\urladdr{http://mmf.ruc.dk/$\sim$booss}

\author{Matthias Lesch}
\address{Mathematisches Institut\\
Universit\"at zu K\"oln\\
Weyertal 86--90\\
D--50931 K\"oln\\
Germany}
\email{lesch@mi.uni-koeln.de}
\urladdr{http://www.mi.uni-koeln.de/$\sim$lesch}
%
\author{John Phillips}
\address{Department of Mathematics and Statistics\\
University of Victoria\\
Victoria, B.C.\\
Canada V8W 3P4}
\email{phillips@math.uvic.ca}


\subjclass[2000]{Primary: 58J30; Secondary: 47A53, 19K56, 58J32
\\
Submitted September 2001, in revised form February 2004} 

\begin{abstract}
We study the gap (= ``projection norm" = ``graph distance") 
topology of the space of all (not necessarily bounded)
self-adjoint Fredholm operators in a separable
Hilbert space by the Cayley transform and direct methods.
In particular, we show the surprising result that this space is
connected in contrast to the bounded case. Moreover, we
present a rigorous definition
of spectral flow of a path of such operators (actually alternative but
mutually equivalent definitions) and prove the homotopy invariance. As an
example, we discuss operator curves on manifolds with boundary. 
\end{abstract}

\maketitle

\section*{Introduction}\label{sec1}
The main purpose of this paper is to study the topology of the space of all 
(generally unbounded) self-adjoint Fredholm operators, and to put the notion 
of spectral flow for continuous paths of such operators
on a firm mathematical footing with clear concise
definitions and proofs. 

The natural topology on the space of all such operators, denoted
by $\cC\!\cF^{\sa}$\,, (for a fixed separable Hilbert space, $H$) is
given by the graph distance topology. That is, we consider the topology 
induced by the metric: $\gd(T_1,T_2)=\|P_1-P_2\|$ where $P_i$ is the
projection onto the graph of $T_i$ in the space $H\times H$ for
$i=1,2$. This metric is called the gap metric. The space of unbounded Fredholm
operators has been studied systematically in the seminal paper by Cordes
and Labrousse \cite{CorLab:IIM}.

Many users of the notion of spectral flow feel that 
the definition and basic properties are already well-understood.
However, there are some difficulties with the currently available
definitions which this paper aims to remedy.

A feature of our approach is the use of the Cayley Transform, 
\[
T\mapsto \Kappa(T) = (T-i)(T+i)\ii\,.
\]
We show that the image $\Kappa(\cC\!\cF^{\sa})$ is precisely the set
\[
\{U\in \cU(H)\mid (U+I) \text{ is Fredholm and $(U-I)$ is
injective}\}=:{_{\cF}\!\cU}_{\rm inj}\,, 
\]
and that the map $\Kappa$
induces an equivalent metric, $\wt\gd$, on $\cC\!\cF^{\sa}$ via
\[
\wt\gd(T_1,T_2) = \|\Kappa(T_1)-\Kappa(T_2)\|.
\]
This Cayley picture of $\cC\!\cF^{\sa}$
leads us to a more careful study of the metric space
$\cC\!\cF^{\sa}$ by studying its image ${_{\cF}\!\cU}_{\rm inj}
 =\Kappa(\cC\!\cF^{\sa})$. In contrast to the space of bounded
self-adjoint Fredholm operators, we  
prove the surprising result that $\cC\!\cF^{\sa}$ is (path-)connected.
In particular, the operator $I$ can be connected to $-I$ in $\cC\!\cF^{\sa}$\,.

Furthermore, using the Cayley picture of $\cC\!\cF^{\sa}$\,, we are able
to give two different (but equivalent) definitions of the spectral
flow of a continuous path in $\cC\!\cF^{\sa}$ and to show that these
definitions are invariant under homotopy. We use neither
Kato's Selection Theorem nor any differentiability or 
regularity assumptions. Thus, spectral flow
induces a surjective homomorphism $\SF$, from the fundamental
group $\pi_1(\cC\!\cF^{\sa})$ to $\Z$.

On the other hand, the space $\cF^{\sa}$  of bounded operators
in $\cC\!\cF^{\sa}$ inherits its usual (norm) topology with the
gap metric $\gd$ and $\cF^{\sa}$ has three connected components by a result
of Atiyah and Singer. To add to the confusion, $\cF^{\sa}$ is also dense
in $\cC\!\cF^{\sa}$. Unfortunately, we have been unable to decide 
whether $\SF:\pi_1(\cC\!\cF^{\sa})\to\Z$ is injective or whether $\cC\!\cF^{\sa}$ 
is a classifying space for the $K^1$-functor 
(cf. Remark \plref{r:connected} below).

Finally, we consider a fixed compact Riemannian manifold $M$ with boundary $\Sigma$, 
a family $\{D_s\}$ of symmetric elliptic differential operators
of first order and of Dirac type on $M$ acting on sections of a
fixed Hermitian vector bundle $E$ with coefficients depending 
continuously on a parameter $s$, and a norm-continuous family $\{P_t\}$ of 
orthogonal projections of $L^2(\Sigma;E|_\Sigma)$ defining well-posed
boundary problems. Here ``Dirac type" means that each operator $D_s$ can be
written in product form near any hypersurface (for details of the 
definition see Assumption \plref{a:dirac} (1), 
Equation \eqref{e:product} below). 

With a view to applications in low-dimensional topology and gauge theories
(see e.g. \cite{KirLes:EIM}), we do not assume that the 
metric structures of $M$ and $E$ are of product form near $\Sigma$.
Consequently the tangential symmetric and skew-symmetric operator
components may depend on the normal variable near $\Sigma$. Solely 
exploiting elliptic regularity and the {\em unique continuation property} of
operators of Dirac type, 
we show that the induced two-parameter family 
\[
(s,t) \mapsto (D_s)_{P_t}
\]
of self-adjoint $L^2(M;E)$-extensions with compact resolvent
is continuous in $\cC\!\cF^{\sa}(L^2(M;E))$ in the gap metric without any 
further assumptions or restrictions.

\medskip 
The results of this paper have been announced in \cite{BooLesPhi:SSF}.

\subsection{Notations}
 Let $H$ be a separable complex Hilbert space. First
let us introduce some notation for various spaces of operators in $H$:
\begin{align*}
 \cC(H) & \text{ closed densely defined operators in $H$},\\
 \cB(H) & \text{ bounded linear operators $H\to H$},\\
 \cU(H) & \text{ unitary operators $H\to H$},\\
 \cK(H) & \text{ compact linear operators $H\to H$},\\
 \cF(H) & \text{ bounded Fredholm operators $H\to H$},\\
 \cC\!\cF(H) & \text{ closed densely defined Fredholm operators in $H$}.
\end{align*}

If no confusion is possible we will omit ``$(H)$'' and write $\cC, \cB,
\cK$ etc. By $\cC^{\sa}, \cB^{\sa}$ etc. we denote the set of self-adjoint
elements in $\cC, \cB$ etc.

\section{The space of unbounded self-adjoint Fredholm operators}

\subsection{The topology of $\cC^{\sa}(H)$} We present a few facts about the
so called gap topology on $\cC^{\sa}$, cf. \cite{CorLab:IIM}, \cite{Kat:PTL},
\cite{Nic:SFO}.
As explained, e.g., in \cite[Section 1]{Nic:SFO} there are two natural metrics
on $\cC^{\sa}$, the {Riesz metric} and the {gap metric}. The \emph{Riesz
metric} is the metric such that the bijection
\begin{equation}\label{G1.1}
\begin{matrix}
   F:& \cC^{\sa} & \longrightarrow & \bigsetdef{S\in\cB^{\sa}}
{\|S\|\le 1 \text{ and $S\pm I$ both injective}},\\ 
\ & T&\mapsto &F_T:=T(I+T^2)^{-1/2}
\end{matrix}
\end{equation}
is an isometry. That is, given $T_1, T_2\in\cC^{\sa}$ then
their Riesz distance $\varphi(T_1,T_2)$ is defined
to be $\|F_{T_1}-F_{T_2}\|$. Note that the image of $F$ is neither
open nor closed in the closed unit ball of $\cB^{\sa}$\,.
Note also that $F$ maps the space $\cC\!\cF^{\sa}$ of (generally unbounded)
self-adjoint Fredholm operators onto the intersection of the space $\cF^{\sa}$ 
of bounded self-adjoint Fredholm operators with $F(\cC^{\sa})$\,, see
also Subsection \plref{ss:character}. It is clear that $F$ is injective.
We postpone the proof that  $F$ as defined in \eqref{G1.1} is surjective
(see Proposition \plref{p:riesz} below).

The \emph{gap metric}
$\delta(T_1,T_2)$ is given as follows: let $P_j$ denote the
orthogonal projections onto the graphs of $T_j$ in $H\times H$. Then
$\delta(T_1,T_2):=\|P_1-P_2\|$. It is shown in \cite[Section 1]{Nic:SFO} that
the Riesz topology is finer than the gap topology. By an example due to
Fuglede (presented in loc. cit.; see also Example \plref{x:fuglede-curve}
below), the Riesz topology is not equal to
the gap topology and hence the Riesz
topology is strictly finer than the gap topology. This means in particular
that the Riesz map $F$ is not continuous on $(\cC^{\sa},\delta)$. This was
also noted in \cite[Section 4.2]{BooFur:MIF}.

The next result shows that, as for the Riesz topology,
the gap topology can be obtained from a map into the bounded linear
operators. 

Recall that two metrics for the same set are {\em (topologically)
equivalent} if and only if they define the same topology and {\em (uniformly)
equivalent} if and only if they can be estimated mutually in a uniform way. In the
latter case the maps ${\rm id}:(X,\delta_1)\to (X,\delta_2)$ and 
${\rm id}:(X,\delta_2)\to (X,\delta_1)$ are Lipschitz continuous and
thus uniformly continuous.

\begin{theorem}\label{S1.1} 
(a) On $\cC^{\sa}$ the gap metric is (uniformly) equivalent to the metric
$\gamma$ given by
\[ \gamma(T_1,T_2)=\|(T_1+i)^{-1}-(T_2+i)^{-1}\|.\]

\noi (b)
Let $\kappa:\R\to S^1\setminus\{1\}, x\mapsto \frac{x-i}{x+i}$
denote the Cayley transform. Then $\kappa$ induces a homeomorphism
\begin{equation}\label{G1.2}\begin{split}
   \Kappa:&
 \cC^{\sa}(H)\longrightarrow 
\bigsetdef{U\in \cU(H)}{U-I \text{ is injective }},\\      
&T\mapsto \Kappa(T)=(T-i)(T+i)^{-1}.
\end{split}
\end{equation}
More precisely, the gap metric is (uniformly) equivalent to the metric
$\tilde \delta$ defined by $\tilde\delta(T_1,T_2)=
\|\Kappa(T_1)-\Kappa(T_2)\|=\frac{1}{2}\gamma(T_1,T_2)$.
\end{theorem}

\begin{proof} First we recall that for $T\in\cC^{\sa}$ the orthogonal
projection $P_T$ onto the graph of $T$ is given by
\begin{equation}\label{G1.3}
\begin{pmatrix} R_T& TR_T\\ TR_T & T^2R_T\end{pmatrix},\quad
R_T:=(I+T^2)^{-1}. \end{equation}
Hence, the gap metric $\delta$ is (uniformly) equivalent to
\begin{equation}\label{G1.4}
   \delta_1(T_1,T_2)=\|R_{T_1}-R_{T_2}\|+\|T_1R_{T_1}-T_2R_{T_2}\|,
\end{equation}
(see also \cite[Lemma 3.10]{CorLab:IIM}).
The identities
\begin{align*}
(T-i)^{-1}&=(T+i)(T^2+I)^{-1}=TR_T+iR_T\,,\\
(T+i)^{-1}&=(T-i)(T^2+I)^{-1}=TR_T-iR_T 
\end{align*}
yield
\begin{equation}\label{G1.5}\begin{split}
  R_T&=\frac{1}{2i}\bigl((T-i)^{-1}-(T+i)^{-1}\bigr),\\
  TR_T&=\frac 12\bigl((T-i)^{-1}+(T+i)^{-1}\bigr),
\end{split}
\qquad T\in \cC^{\sa},
\end{equation}
and we infer that the metric $\delta_1$ is (uniformly) equivalent to the
metric $\gamma$ given by
\begin{equation}\label{G1.6}\begin{split}
  \gamma(T_1,T_2)&=\frac
12\bigl(\|(T_1+i)^{-1}-(T_2+i)^{-1}\|+\|(T_1-i)^{-1}-(T_2-i)^{-1}\|\bigr)\\ 
    &=\|(T_1+i)^{-1}-(T_2+i)^{-1}\|.
\end{split}
\end{equation}
In the last equality we have used that for any $A\in\cB(H)$ one has
$\|A\|=\|A^*\|$. This proves (a).

\smallskip
To prove (b) we note for $T\in\cC^{\sa}$ the identities $\mathrm{range}(T+i)=H$ and
\begin{equation}\label{G1.7} \Kappa(T)=I-2i (T+i)^{-1}\,.
\end{equation}
These imply
\begin{equation}\label{G1.8}
    \|(T_1+i)^{-1}-(T_2+i)^{-1}\|=\frac 12\|\Kappa(T_1)-\Kappa(T_2)\|.
\end{equation}
This shows that the gap metric and the metric $\tilde \delta$ are
(uniformly) equivalent. This equivalence implies that the Cayley transform
is a homeomorphism onto its image. It remains to identify the image of the
Cayley transform.

Given $T\in\cC^{\sa}$ its Cayley transform $\Kappa(T)$ is certainly a unitary
operator. To show that $\Kappa(T)-I$ is injective consider $x\in H$ such
that $\Kappa(T)x=x$. In view of \eqref{G1.7} this implies
\[x=\Kappa(T)x=x-2i (T+i)^{-1}x;\]
thus $(T+i)^{-1}x=0$ and hence $x=0$.

Conversely, let $U$ be a unitary operator such that $U-I$ is injective.
From the following proposition and corollary, we obtain the existence of a
$T\in\cC^{\sa}$  such that $\Kappa(T)=U$. 
The theorem is proved.\end{proof}

\begin{prop}\label{p:cayley-inv}
If $U$ is unitary and $U-I$ is injective, then $T:= i(I+U)(I-U)\ii$ is
self-adjoint on $\dom(T):=\mathrm{range}(I-U)$. Moreover, $T=i(I-U)\ii(I+U)$.
\end{prop}

A similar result is proved in \cite[Theorem 13.19]{Rud:FA}. Our
argument seems to be shorter and more appropriate
in our context.

\begin{proof}
$\overline{\mathrm{range}(I-U)} = \ker (I-U^*)^\perp = \ker (I-U)^\perp = \{0\}^\perp 
= H$ since $U$ normal implies $\ker(I-U^*) = \ker (I-U)$. 
Thus, $\dom T$ is dense in $H$.
Now, 
\begin{align*}
(I+U)(I-U)\ii &= (I-U)\ii(I-U)(I+U)(I-U)\ii \\
&= (I-U)\ii (I+U)|_{\mathrm{range}(I-U)} \subseteq (I-U)\ii(I+U).
\end{align*}
On the other hand, if $x\in\dom\bigl((I-U)\ii(I+U)\bigr)$ then
\[
(I+U)x\in\dom\bigl((I-U)\ii\bigr) = \mathrm{range}(I-U),
\]
and so there exists a $y\in H$
with $(I+U)x=(I-U)y$. Solving, 
\[
x=(I-U)y + (I-U)x - x
\]
and so $x=(I-U)\12(x+y)\in\dom\bigl[(I+U)(I-U)\ii\bigr]$. Thus,
\[
T= i(I+U)(I-U)\ii = i(I-U)\ii(I+U).
\]

It is an elementary calculation that $T$ is symmetric and so 
\[
T \subseteq T^* = -i(I-U^*)\ii(I+U^*)
\]
(we have the ``=" 
since $I+U$ is bounded and on the left in the
formula for $T$, see e.g. \cite[p. 299]{RieNag:FA}) and by the
same argument as for $T$ we get
\[
T^* = -i(I-U^*)\ii(I+U^*) = -i(I+U^*)(I-U^*)\ii
\]
and $T^*$ is symmetric, so that
\[
T^* \subseteq T^{**} = i(I-U)\ii(I+U) = T.
\]
Hence, $T=T^*$\,.
\end{proof}

\begin{cor}
With $U$ and $T$ as above, $\Kappa(T)=U$.
\end{cor}

\begin{proof}
\begin{align*}
(T+iI) &= i(I-U)\ii(I+U) + i(I-U)\ii (I-U)\\
&= i(I-U)\ii \cdot 2 = 2i(I-U)\ii,
\end{align*}
so that 
\[
(T+iI)\ii = \frac 1{2i}(I-U).
\]
By a similar calculation,
\[
(T-iI)=2i(I-U)\ii U=2iU(I-U)\ii
\]
so that,
\[
\Kappa(T)=(T-iI)(T+iI)\ii=U.
\]
\end{proof}

\begin{remark}\label{r:metrics}
(a) In the definition of the metric $\gamma$ in
\eqref{G1.6} we may replace $i$ by $-i$ or, more generally, 
by any $-\lambda$ with $\lambda\in {\varrho}(T_1) \cap {\varrho}(T_2)$, 
$\varrho(T_j):=\C\setminus \spec T_j$ denoting the resolvent set. 
All these metrics are (uniformly) equivalent with the gap metric.

\noi (b) We recall the basic spectral argument for Cayley transforms,
namely that the identity $\gl I-T = (\gl+i) (\kappa(\gl)-\Kappa(T))
(I-\Kappa(T))\ii$ implies
\begin{align}
\gl\in\spec T \quad &\iff\quad \kappa(\gl)\in\spec\Kappa(T),\\
\gl \in \spec_{{\rm discr}} T  \quad &\iff\quad 
 \kappa(\gl) \in\spec_{{\rm discr}} 
\Kappa(T)\,.
\end{align} 
Here $\spec_{{\rm disc}}$ denotes the discrete spectrum, cf. 
subsection \plref{ss:character} below.
\end{remark}

\medskip

Following the same pattern as the preceding proof of Proposition 
\plref{p:cayley-inv} we show

\begin{prop}\label{p:riesz}
If $S$ is a bounded self-adjoint operator with $\|S\|\le 1$ and $S\pm I$ 
injective, then $T:=S(I-S^2)^{-\12}$ is densely defined and self-adjoint.
Moreover, 
\[
T=(I-S^2)^{-\12}\/S \ \tand\ S=T(I+T^2)^{-\12}\,.
\]
\end{prop}

\begin{proof} Since $I-S^2$ is injective it has dense range and so
$(I-S^2)\ii$ and $(I-S^2)^{-\12}$ are densely defined
and self-adjoint. Since $S$ commutes with $(I-
S^2)^{\12}$ we have that $S(I-S^2)^{-\12}\subseteq
(I-S^2)^{-\12}\/S$ by an argument in the proof of Proposition 
\plref{p:cayley-inv}. On the other hand, for
$x\in\dom \bigl((I-S^2)^{-\12}\/ S\bigr)$ we have
$Sx\in \dom \bigl((I-S^2)^{-\12} \bigr) =\mathrm{range} \bigl((I-
S^2)^{\12} \bigr)$ so that 
\[
Sx=(I-S^2)^{\12}y 
\]
for some $y$. Hence, $S^2x= S(I-S^2)^{\12}y = (I-S^2)^{\12}\/Sy$.
Or, $(I-S^2)x= x-(I-S^2)^{\12}\/Sy$. That is, 
\[
x = (I-S^2)x+(I-S^2)^{\12}\/Sy
= (I-S^2)^{\12}\bigl((I-S^2)^{\12}x+Sy\bigr)
\]
is in the range of $(I-S^2)^{\12}$ which is 
$\dom\bigl((I-S^2)^{-\12}\bigr) = \dom\bigl(S(I-S^2)^{-
\12}\bigr)$. That is, $(I-S^2)^{-\12}\/S=S(I-S^2)^{-\12}$\,.
By an argument in the proof of Proposition 
\plref{p:cayley-inv}, this implies that $T:=(I-S^2)^{-
\12}\/S$ is self-adjoint.

Now, since $S$ commutes with $(I-S^2)^{-\12}$ one
calculates
\[
(I+T^2) = I + (I-S^2)\ii S^2 = (I-S^2)\ii \bigl((I-
S^2)+S^2\bigr) = (I-S^2)\ii \,.
\]
From this we easily calculate $T(I+T^2)^{-\12} = S$.
\end{proof}

\medskip

It was proved in \cite[Addendum]{CorLab:IIM} that the topology induced by
the gap metric on the set of bounded operators is the same as the topology
induced by the natural metric $s(T_1,T_2)=\|T_1-T_2\|$. However, the reader
should be warned that the metric $s$ is not (uniformly) equivalent to the
gap metric. In other words, the uniform structures induced by the gap metric
and by the operator norm on the space of bounded linear operators are
different.
This follows from the fact that the metric $s$ is complete while the gap
metric on the set of bounded operators is not complete. The latter
follows from the following result.

\begin{prop}\label{S1.2} With respect to the gap metric
the set $\cB^{\sa}(H)$ is dense
in $\cC^{\sa}(H)$.
\end{prop}
\begin{proof} Let $T\in\cC^{\sa}$ and denote by
$(E_\lambda)_{\gl\in\R}$ the spectral
resolution of $T$. Put
\begin{equation}
  T_n:=\int_{[-n,n]}\gl dE_\gl+\int_{|\lambda|>n} 
    n (\sgn\gl)dE_\gl\,. 
\end{equation}
Then $T_n$ is a bounded self-adjoint operator and
\begin{equation}\label{G1.10}
    \begin{split}
    \gamma(T,T_n)&=
\|(T+i)\ii-(T_n+i)\ii\|\\
      &=\|\int_{|\gl|>n}(\gl+i)^{-1}-(n(\sgn\gl)+i)^{-1}dE_\gl\|      \le \frac
2n\,.   
   \end{split}
\end{equation}
Hence $T_n\to T$ in the $\gamma$-metric.
In view of Theorem \ref{S1.1} (a) this proves the assertion.
\end{proof}

\subsection{The connectedness of $\cC\!\cF^{\sa}$}\label{ss:character}
We determine the image under the Cayley transform of the space $\cC\!\cF^{\sa}$ of (not
necessarily bounded) self-adjoint Fredholm operators. Moreover, we will
show that this space is path connected. For the general
theory of
unbounded Fredholm operators we refer to \cite[Section IV.5]{Kat:PTL}.

We recall that for a closed operator $T$ in a Hilbert
space the essential spectrum, $\spec_{{\rm ess}} T$,
consists of those $\lambda\in\C$ for which $T-\lambda$ is not a Fredholm
operator. Then $\spec_{{\rm ess}} T$ is a closed subset of $\spec T$. The
discrete spectrum, $\spec_{\rm discr} T$, consists of those isolated
points of $\spec T$
which are not in $\spec_{\rm ess} T$.

It is well-known that if $T$ is self-adjoint then
 $\lambda$ is an isolated point of $\spec T$ if and only if 
$\mathrm{range} (T-\lambda)$ is
closed (\cite[Definition XIII.6.1 and Theorem XIII.6.5]{DunSchwa:LO}; 
note that loc. cit.
define the essential spectrum differently). Consequently, for a
self-adjoint operator $T$ we have 
\[
\begin{split} 
\spec_{\rm discr} T&=\spec T\setminus \spec_{\rm ess} T\\
       &=\{\lambda\in\C\,|\, \lambda \text{ is an isolated point of }
\spec T \text{ which is}\\
&\qquad\qquad\qquad\text{an eigenvalue of {\em finite} multiplicity of } T\}\\
       &=\{\lambda\in\C\,|\, 0 < \dim \ker (T-\lambda)<\infty \text{ and 
$\mathrm{range} (T-\lambda)$  closed }\}.
\end{split} 
\]

\medskip

We note an immediate consequence of the Cayley picture:

\begin{prop}\label{S2.1} For $\gl\in\R$ the sets
\[ \bigsetdef{T\in\cC^{\sa}(H)}{\gl\not\in\spec T} \qquad \text{and} \qquad   
\bigsetdef{T\in\cC^{\sa}(H)}{\gl\not\in\specess T}
\]
are open in the gap topology.
\end{prop}

\begin{proof} By Theorem \ref{S1.1} (see also Remark \ref{r:metrics}b) we
have
\[
\begin{split}
 \bigsetdef{T\in\cC^{\sa}(H)}{\gl\not\in\spec T}
&=\Kappa^{-1}\bigsetdef{U\in\cU(H)}{\kappa(\gl)\not\in\spec U},\\   
\bigsetdef{T\in\cC^{\sa}(H)}{\gl\not\in\specess T}&=
     \Kappa^{-1}\bigsetdef{U\in\cU(H)}{\kappa(\gl)\not\in\specess U},
  \end{split}
\]
where the spaces of unitary operators on the right side are open in the
range of $\Kappa$ by the openness of the spaces of bounded invertible resp.
bounded Fredholm operators. Now the assertion follows.
\end{proof}

\begin{cor}\label{c:fred-open}
The set $\cC\!\cF^{\sa}=\bigsetdef{T\in\cC^{\sa}}{0\not\in\specess T}
=\Kappa^{-1}({{_\cF}\cU})$, ${{_\cF}\cU}:=\bigsetdef{U\in\cU}{-1\not\in\specess
U}$, of (not necessarily bounded) self-adjoint Fredholm operators is open
in $\cC^{\sa}$.
\end{cor}

\begin{remark}\label{r:open}
By Proposition \plref{S1.2}, the preceding corollary implies that the set
$\cF^{\sa}$ is dense in $\cC\!\cF^{\sa}$ with respect to the gap metric.
\end{remark}

Contrary to the bounded case and somewhat surprisingly
the space of unbounded self-adjoint Fredholm operators is connected. 
More precisely we have:

\begin{theorem}\label{t:connected}
(a) 
The set $\cC\!\cF^{\sa}$ is path-connected with respect to the gap metric.

\noi (b) Moreover, its Cayley image 
\[
{_{\cF}\!\cU}_{\!\inj} := \bigsetdef{U\in\cU}{U+I \text{ Fredholm and $U-I$
injective}} = \Kappa(\cC\!\cF^{\sa}) 
\]
is dense in ${_{\cF}\!\cU}$.
\end{theorem}

\begin{proof}
(a) Once again we look at the Cayley transform picture. We shall use the
following notation:
\[
\cU_{\!\inj} := \bigsetdef{U\in\cU}{U-I \text{ injective}}
= \Kappa(\cC^{\sa}). 
\]
Note that ${_{\cF}\!\cU}_{\!\inj} = {_{\cF}\!\cU}\cap \cU_{\!\inj}$\,. 
We consider a fixed $U\in {_{\cF}\!\cU}_{\!\inj}$\,. Then $H$ is the direct sum
of the spectral subspaces $H_\pm$
of $U$ corresponding to $[0,\pi)$ and $[\pi,2\pi]$ respectively and we may 
decompose $U=U_+\oplus U_-$. More precisely, we have
\[
\spec(U_+) \< \bigsetdef{e^{it}}{t\in [0,\pi)}
\text{ and }
\spec(U_-) \< \bigsetdef{e^{it}}{t\in [\pi,2\pi]}\,.
\]
Note that there is no intersection between the spectral spaces in the
endpoints: if $-1$ belongs to $\spec(U)$, it is an isolated eigenvalue by
our assumption and hence belongs only to $\spec(U_-)$; if $1$ belongs to 
$\spec(U)$, it can belong both to $\spec(U_+)$ and $\spec(U_-)$, 
but in any case, it does not
contribute to the decomposition of $U$ since, by our assumption, $1$ is not
an eigenvalue at all. 

By spectral deformation (squeezing the spectrum down to
$+i$ and $-i$) we contract $U_+$ to $iI_+$ and $U_-$ to $-iI_-$\,, where
$I_\pm$ denotes the identity on $H_\pm$\,. We do this on the upper half arc
and the lower half arc, respectively, in such a way that 1 does not become
an eigenvalue under the course of the deformation: actually it will no
longer belong to the spectrum; neither will $-1$ belong to the spectrum.
That is, we have connected $U$ and $iI_+\oplus -iI_-$ within
$\Kappa(\cC\!\cF^{\sa})$.

\begin{figure}
\input{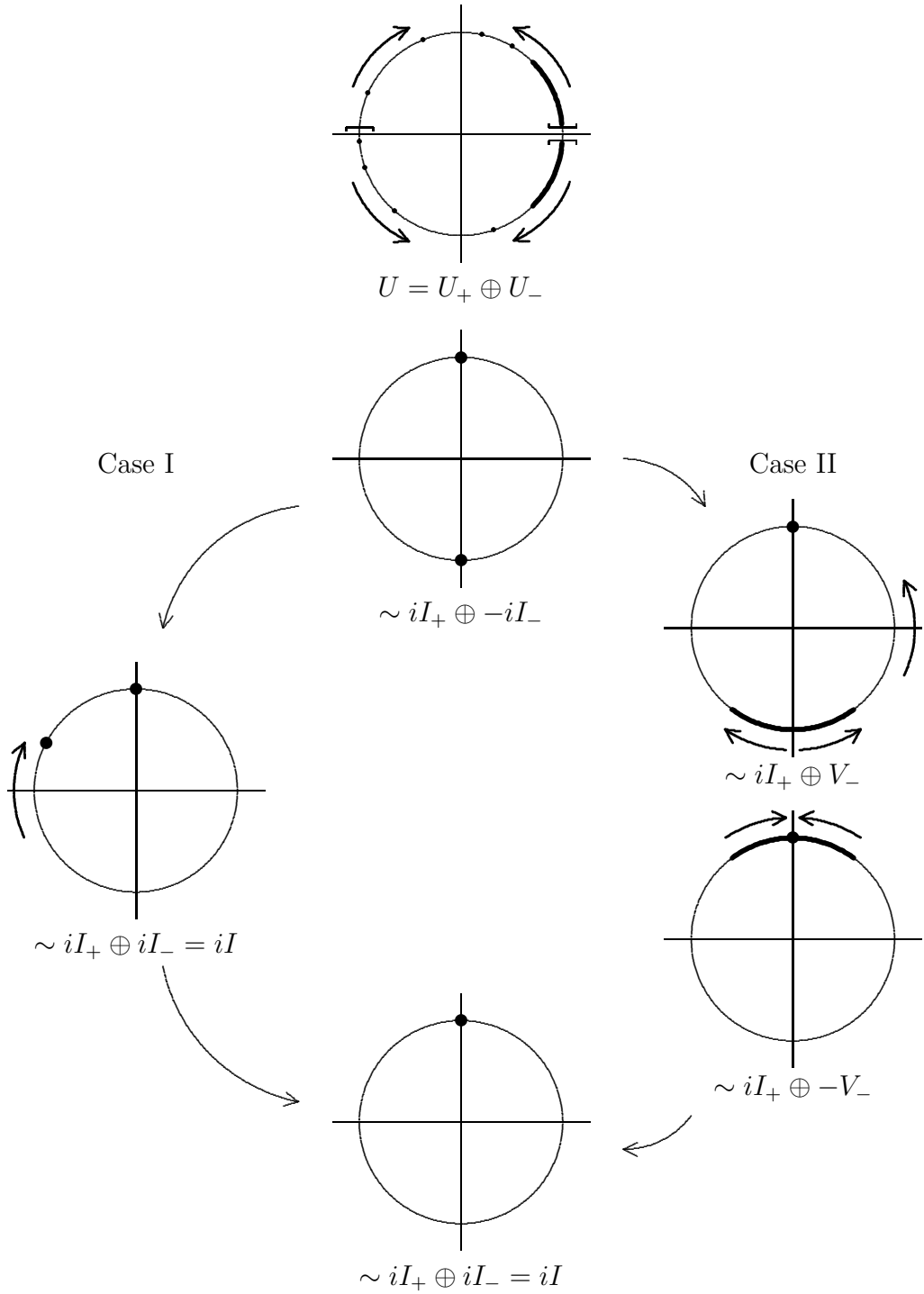}
\caption{\label{f:deformation} Connecting a fixed $U$ in
${_{\cF}\!\cU}_{\!\inj}$ to $iI$. Case I (finite rank $U_-$) and
Case II (infinite rank $U_-$)}
\end{figure}

We distinguish two cases:
If $H_-$ is finite-dimensional, we now rotate $-iI_-$ up through $-1$
into $iI_-$\.. More precisely, we consider $\{iI_+ \oplus e^{i(\pi/2
+(1-t)\pi)}I_-\}_{t\in [0,1]}$\,. This proves that we can connect $U$ with
$iI_+\oplus iI_- = iI$ within $\Kappa(\cC\!\cF^{\sa})$ in this first case.

If $H_-$ is infinite-dimensional, we ``un-contract" $-iI_-$ in such
a way that no eigenvalues remain. To do this, we identify $H_-$ with
$L^2([0,1])$. Now multiplication by $-i$ on $L^2([0,1])$
can be connected to multiplication by a function whose values are
a short arc centred on $-i$ and so that the resulting operator $V_-$ on
$H_{-}$    has no eigenvalues. This will at no time introduce spectrum near
$+1$ or $-1$. We then rotate this arc up through $+1$ (which keeps us in
the right space) until it is centred on $+i$. Then we contract the spectrum
on $H_{-}$ to be $+i$. That is, also in this case we have connected our
original operator $U$ to $+iI$. To sum up this second case 
(see also Figure \ref{f:deformation}):
\begin{multline*}
U  \sim iI_+ \oplus -iI_- \sim iI_+ \oplus V_-
\sim iI_+ \oplus e^{it\pi}V_- \text{ for $t\in [0,1]$}\\
\sim iI_+ \oplus -V_- \sim iI_+ \oplus -(-iI_-) \sim iI\,.
\end{multline*} 

\noi To prove (b), we just decompose any $V\in {_{\cF}\!\cU}$ into $V=U\oplus
I_1$ where $U\in {_{\cF}\!\cU}_{\!\inj}(H_0)$ and $I_1$ denotes the identity on
the 1-eigenspace $H_1=\ker (V-I)$ of $V$ with $H=H_0\oplus H_1$ an
orthogonal decomposition. Then  for $\eps>0$, $U\oplus e^{i\eps} I_1 \in
{_{\cF}\!\cU}_{\!\inj}$ approaches $U$ for $\eps\to 0$.
\end{proof}

\begin{remark}\label{r:connected} 
Recall that $\cF^{\sa}$ has three connected components
\[
  \cF^{\sa}_\pm=\bigsetdef{T\in\cF^{\sa}}{\specess(T)\subset \R_\pm},
\]
and $\cF^{\sa}_*=\cF^{\sa}\setminus\Bigl(\cF^{\sa}_+\cup \cF^{\sa}_-\Bigr)$.
$\cF^{\sa}_\pm$ are contractible and $\cF^{\sa}_*$ is a classifying space
for the $K^1$--functor \cite{AtiSin:ITS}. Recall that K--theory is a
generalized cohomology theory; a classifying space for $K^1$ is
a topological space $Z$ such that $K^1(X)$ is naturally isomorphic
to the homotopy classes, $[X,Z]$, of maps $X\to Z$.

The preceding proof shows also that the two subsets of $\cC\!\cF^{\sa}$
\[
\cC\!\cF_\pm^{\sa} = \bigsetdef{T\in\cC\!\cF^{\sa}} 
{\specess(T)\< \R_{\pm}},
\]
the spaces of all essentially positive resp. all essentially
negative self-adjoint Fredholm operators, are
no longer open. The third of the three complementary subsets 
\begin{equation}
\cC\!\cF_*^{\sa} = \cC\!\cF^{\sa} \setminus 
\Bigl(\cC\!\cF_+^{\sa} \cup \cC\!\cF_-^{\sa} \Bigr)
\label{revision:1.13}
\end{equation}
is also not open. 
We do not know whether the two ``trivial" components are
contractible as in the bounded case nor whether the whole space is a
classifying space for $K^1$ as is the non-trivial component in the bounded
case.


Independently of the Fuglede example, the connectedness of $\cC\!\cF^{\sa}$ and
the non-connectedness of $\cF^{\sa}$ show 
that the Riesz map is not continuous on $\cC\!\cF^{\sa}$ in the gap topology.
\end{remark}

\section{Spectral flow for unbounded self-adjoint operators}

\subsection{First approach via Cayley transform and winding number}

In \cite[Section 6]{KirLes:EIM} it was shown that
the natural inclusion 
\[
\cU_\cK(H):=\bigsetdef{U\in\cU}{U-I \text{ is
compact }} \hookrightarrow
{_{\cF}\!\cU}(H):=\bigsetdef{U\in\cU}{-1\not\in\specess U}
\]
is a homotopy
equivalence. As a consequence the classical winding number extends to an
isomorphism
\begin{equation}
  \wind: \pi_1({_{\cF}\!\cU},I)\longrightarrow \Z\,,
\end{equation}
see also \cite[Appendix]{FurOts:MII} for a different proof (cf. also
Proposition \plref{p:retract} below).

Furthermore, in \cite[l.c.]{KirLes:EIM}  it was shown that to any continuous (not
necessarily closed) curve $f:[0,1]\to {_{\cF}\!\cU}$ one can assign an
integer $\wind(f)$ in such a way that the mapping $\wind$ is
\begin{enumerate}
\item \textit{Path additive}: Let $f_1, f_2:[0,1]\to{_{\cF}\!\cU}(H)$
be continuous paths with 
\[f_2(0)=f_1(1).\] 
Then
\[\wind(f_1*f_2)=\wind(f_1)+\wind(f_2).\]
\item \textit{Homotopy invariant}: Let $f_1,f_2$ be continuous paths
in ${_{\cF}\!\cU}$. Assume that there is a homotopy
$H:[0,1]\times[0,1]\to{_{\cF}\!\cU}$
such that $H(0,t)=f_1(t), H(1,t)=f_2(t)$ and such that
$\dim \ker(H(s,0)+I), \dim \ker(H(s,1)+I)$ are independent of $s$.
Then
$\wind(f_1)=\wind(f_2)$. In particular, $\wind$ is invariant under
homotopies leaving the endpoints fixed.
\end{enumerate}
Roughly speaking, the mapping $\wind$ is the `spectral flow' across $-1$; 
that is, $\wind$
counts the net number of eigenvalues of $f(t)$ which cross $-1$ from the
upper half-plane into the lower half-plane. One has to choose a convention
for those cases in which $-1\in\spec f(0)$ or $-1\in\spec f(1)$. Contrary to 
the convention which was chosen in \cite{KirLes:EIM}, our convention is chosen 
as follows: choose $\eps>0$ so small that
$-1\not\in\spec(f(j)e^{i\varphi}), j=0,1$ for all $0<|\varphi|\le\eps$.
Then put $\wind(f):=\wind(fe^{i\eps})$. This means that an eigenvalue
running from the lower half-plane into $-1$ is not counted while an
eigenvalue running from the upper half-plane into $-1$ contributes $1$ to
the winding number.

In analogy to \cite{Phi:SAF} we can give an explicit description of
$\wind(f)$. Alternatively, it can be used as a definition of $\wind$:

\begin{prop}\label{p:wind-def}
Let $f:[0,1]\to{_{\cF}\!\cU}$ be a continuous path. 

\noi (a) There is a partition $\{0=t_0<t_1<\dots<t_n=1\}$ of the interval
and positive real numbers $0<\eps_j<\pi$, $j=1,\dots, n$, such that
$\ker(f(t)-e^{i(\pi\pm\eps_j)})=\{0\}$ for $t_{j-1}\le t\le t_j$\,. 

\noi (b) Then
\begin{equation}\label{e:wind-def}
\wind(f)=\sum_{j=1}^n k(t_j,\eps_j)-k(t_{j-1},\eps_j),
\end{equation}
where 
\[
k(t,\eps_j):= \sum_{0\le\theta<\eps_j} \dim\ker(f(t)-e^{i(\pi+\theta)}). \]

\noi (c) In particular, this calculation of $\wind(f)$ is independent of
the choice of the segmentation of the interval and of the choice of the
barriers.
\end{prop}

\begin{proof}
In (a) we use that $f(t)\in{_{\cF}\!\cU}$ and $f$ continuous. (b) follows
from the path additivity of $\wind$. (c) is immediate from (b).
\end{proof}

This idea of a {\em spectral flow across $-1$} was introduced first
in \cite[Section 1.3]{BooFur:MIF}, where it was used to give a definition of
the Maslov index in an infinite dimensional context.

After these explanations the definition of spectral flow for paths in
$\cC\!\cF^{\sa}$ is straightforward:

\begin{dfn}\label{S2.2} Let $f:[0,1]\to \cC\!\cF^{\sa}(H)$ be a continuous
path. Then the {\em spectral flow} of $f$, $\SF(f)$ is defined by  
\[
    \SF(f):=\wind(\Kappa\circ f).
\]
\end{dfn}

From the properties of $\Kappa$ and of the winding number we infer
immediately:

\begin{prop}\label{S2.3}
$\SF$ is path additive and homotopy invariant in the following
sense: let $f_1,f_2:[0,1]\to\cC\!\cF^{\sa}$ be continuous paths
and let
\[
H:[0,1]\times[0,1]\to\cC\!\cF^{\sa}
\]
be a homotopy
such that $H(0,t)=f_1(t), H(1,t)=f_2(t)$ and such that
$\dim \ker H(s,0),$ $\dim \ker H(s,1)$ are independent of $s$.
Then $\SF(f_1)=\SF(f_2)$. In particular, $\SF$ is in\-var\-i\-ant under
homotopies leaving the endpoints fixed.
\end{prop}

From Proposition \ref{p:wind-def} we get
\begin{prop} \label{S2.4} For a continuous path $f:[0,1]\to \cF^{\sa}$ our
definition of spectral flow coincides
with the definition in \cite{Phi:SAF}.
\end{prop}

Note that also the conventions coincide for $0\in\spec f(0)$
or $0\in\spec f(1)$.

Returning to the Cayley picture, we have that the mapping $\wind$ induces a surjection
of $\pi_1({_{\cF}\!\cU}_{\!\inj})$ onto $\Z$. Because $\Z$ is free, there is a
right inverse of $\wind$ and a normal subgroup $G$ of
$\pi_1({_{\cF}\!\cU}_{\!\inj})$ such that we have a split short exact sequence
\begin{equation}\label{e:exact}
0\too G \too \pi_1({_{\cF}\!\cU}_{\!\inj}) \too \Z \too 0.
\end{equation}
For now, an open question is whether $G$ is trivial: does the winding
number distinguish the homotopy classes? That is, the question is whether each
loop with winding number 0 can be contracted to a constant point, or,
equivalently, whether two continuous paths in $\cC\!\cF^{\sa}$ with same
endpoints and with same spectral flow can be deformed into each other? Or
is \quad
$\pi_1({_{\cF}\!\cU}_{\!\inj}) \ \cong\ \Z \,\times$\!\!\,{\small\small\small 
$|$}\ $G$\quad the semi-direct product of a non-trivial factor $G$ with $\Z$?

We know a little more than \eqref{e:exact}:

\begin{prop}\label{p:retract}
There exists a continuous map ${_{\cF}\!\cU} \to \cU_\infty$ which induces an
isomorphism $\pi_1({_{\cF}\!\cU}) \to \pi_1(\cU_\infty)=\Z$. Moreover, the
restriction of this map to  ${_{\cF}\!\cU}_{\!\inj}$ induces a map such that
the following diagram commutes
\begin{equation}\label{e:CD_trivialization}
\begin{matrix}
\pi_1({_{\cF}\!\cU}_{\!\inj}) & \quad\too & \pi_1({\cU}_\infty)    \\
\quad &\ &\  \\
\ & \wind\searrow\quad & \quad{\cong}\downarrow{\wind} \\
\quad &\ &\  \\
\ &\ &\Z
\end{matrix}
\end{equation}
\end{prop}

\begin{proof} Let $U_0\in {_{\cF}\!\cU}$\,. Then there exists a neighbourhood
$N_{\eps_0}$ of $U_0$ in ${_{\cF}\!\cU}$ and $\eps_0>0$ such that for each
$U\in N_{\eps_0}$ the projection $\chi_{\eps_0}(U)$ has finite rank where
$\chi_{\eps_0}$ denotes the characteristic function of the arc 
$\{e^{it}\mid t\in[\pi-\eps_0,\pi+\eps_0]\}$ of the unit circle
$\mathbb{T}$. Now, there is a continuous function $f_{\eps_0}: 
\mathbb{T} \to \mathbb{T}$ such that:
\[
f_{\eps_0}(z)=
\begin{cases}
z & \text{ for $z\in\{e^{it}\mid t\in [\pi-\frac{\eps_0}2,
\pi+\frac{\eps_0}2]\}$} \\
1 & \text{ for $z\in
\{e^{it}\mid t\in [0,\pi-\eps_0] \cup [\pi+\eps_0,2\pi]\}$} 
\end{cases}
\]
with
\[
f_{\eps_0}: \begin{cases}
\{e^{it}\mid t\in [\pi-\eps_0,\pi-\frac{\eps_0}2]\}
\to \{e^{it}\mid t\in [0,\pi-\eps_0]\}  \text{ is injective}\\
\{e^{it}\mid t\in [\pi+\frac{\eps_0}2,\pi+\eps_0]\}
\to \{e^{it}\mid t\in [\pi+\eps_0,2\pi]\}  \text{ is injective}.
\end{cases}
\]
Then, actually, $U\mapsto f_{\eps_0}(U)$ : $N_{\eps_0}\to\cU_\infty$\,!

Since ${_{\cF}\!\cU}$ is metric it is paracompact and so the open cover
$\{N_{\eps_0}(U)\}$ has an open locally finite refinement, say $\{N_\ga\}$
and each $N_\ga$ carries a function $f_\ga:N_\ga\to\cU_\infty$ given by a
function $f_\ga: \mathbb{T} \to \mathbb{T}$ corresponding to a positive
$\eps_0$\,. We let $\{p_\ga\}$ be a partition of unity subordinate to the
cover. Then $f: {_{\cF}\!\cU} \to\cB(H)$ is continuous where $f(U):=\sum_\ga
p_\ga(U)f_\ga(U)$. We claim that $f(U)$ is normal and invertible so that
$g(U)=f(U)|f(U)|\ii$ is unitary. To see this, we observe that for each
single $U$ we have $f(U)= \sum_{i=1}^n \gl_i 
f_{\ga_i}(U)$ with the $f_{\ga_i}$ as
above. Moreover, if we let $\gd$ denote the minimum of the corresponding
$\{\12\eps_{\ga_i}\}$ then  $h
=  \sum_{i=1}^n \gl_i f_{\ga_i}$ satisfies 
\begin{align}
h(z)&= z \quad\text{for all $z\in\{e^{it}\mid t\in[\pi-\gd,\pi+\gd]\}$}\\
h(z)&= 1 \quad\text{for all $z\in\{e^{it}\mid t\in [0,\Delta]\cup [2\pi- 
\Delta,2\pi]\}$},
\end{align}
where $\Delta= \max \eps_{\ga_i} >0$ and $\chi_{[\Delta, 2\pi-\Delta]}(U)$
is of finite rank; $h(z)$ lies in one of the shaded convex regions of Figure 
\ref{f:shaded} for all other $z$ on the circle.

\begin{figure}
\input{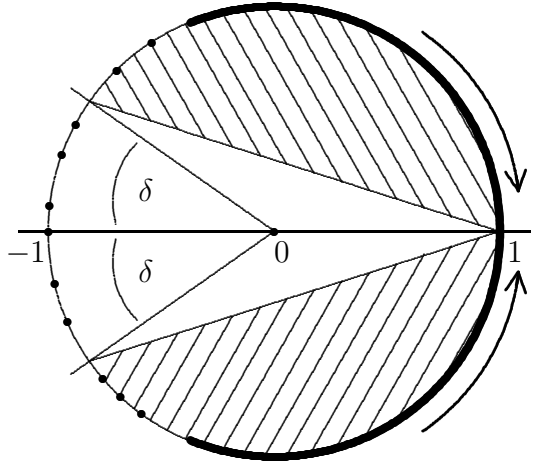}
\caption{\label{f:shaded} Convex regions of finite linear combinations}
\end{figure}

Thus, $f(U)=h(U)$ is normal and invertible. Moreover, since each
\[
f_{\ga_i}(U)\in\cU_\infty \< \{I+ \text{ finite rank operators}\}, 
\]
$f(U)$ is in 
$\{I+ \text{ finite rank operators}\}$ so that $g(U)=f(U) |f(U)|\ii$ is in
$\cU_\infty$\,. Moreover, clearly $\chi_\gd(U) = \chi_\gd(g(U))$ and so we
get the commuting diagram (as the covering is neighbourhood-finite we get
$\chi_\gd(V) = \chi_\gd(g(V))$ for $V$ in a neighbourhood of $U$.
\end{proof}

Summing up it remains an open problem to determine the fundamental group of 
the space ${\mathcal C}{\mathcal F}^{\rm sa}(H)$ or, even more, to determine whether, as
in the bounded case, it is a classifying space for $K^1$. 

Robbin and Salamon \cite{RobSal:SFM} introduced the spectral flow for a family of
unbounded self-adjoint operators under the assumption that the domain is fixed
and that each operator of the family has a compact resolvent. Along the lines
of their method one can prove the following generalization of 
\cite[Theorem 4.25]{RobSal:SFM}:

\begin{prop} Let $f:[0,1]\to {\mathcal C}{\mathcal F}^{\rm sa}(H)$ be a closed 
continuous path. Then
there is a continuous path of self-adjoint matrices 
$g:[0,1]\to {\rm Mat}(n,\C)$ such that 
$f\oplus g$ is homotopic to a closed continuous path of \emph{invertible}
operators 
\[
h:[0,1]\to {\mathcal C}{\mathcal F}^{\rm sa}(H\oplus \C^n).
\]
\end{prop}

If $h$ were a family of bounded invertible operators then it would be clear that it is
homotopic to a constant path. Unfortunately, this is not clear for a path of
unbounded operators. If we could conclude that $h$ is homotopic to a constant 
path
then we would know at least that the ``stable" fundamental group 
of ${\mathcal C}{\mathcal F}^{\sa}(H)$ is isomorphic to $\Z$.


\subsection{Second approach, after \cite{Phi:SAF}}

There is another way of looking at continuous curves of self-adjoint 
Fredholm operators which 
more closely resembles what is done in the bounded self-adjoint
setting. The fact that one can (continuously) isolate the spectra of the
unbounded Fredholm operators in an open interval about 0 is quite appealing
from an operator algebra point of view: it is surprising that this
can be done without the Riesz map being continuous.
Therefore both approaches are included in this note.

In \cite{Phi:SAF} the third author introduced a new method to
define spectral flow of a continuous family of bounded operators. The
interesting new feature of his approach was that it works directly for any
continuous family without first changing the family to a generic situation
(see also Proposition \ref{p:wind-def} above).

In this subsection we adapt the method of \cite{Phi:SAF} to unbounded
operators.

\begin{lemma}\label{S2.5} Let $K\subset\C$ be a compact set. Then
$\bigsetdef{T\in\cC^{\sa}}{K\subset {\varrho}(T)}$ is open in the gap
topology. Here, ${\varrho}(T):=\C\setminus\spec T$ denotes the resolvent
set of $T$.

Similarly, $\bigsetdef{T\in\cC^{\sa}}{K\subset {\varrho}_{\mathrm{ess}}(T)}$,
${\varrho}_\mathrm{ess}(T):=\C\setminus \specess(T)$, 
is open in the gap topology.
\end{lemma}
\begin{proof}
In view of Theorem \ref{S1.1} we find
\begin{equation}\begin{split}
     \bigsetdef{T\in\cC^{\sa}}{K\subset {\varrho}(T)}
    =&\bigsetdef{T\in\cC^{\sa}}{\spec T\subset K^c\cap\R}\\
    =&\bigsetdef{T\in\cC^{\sa}}{\spec\Kappa(T)
\subset\kappa(K^c\cap\R)\cup\{1\}}\\     =&\Kappa^{-1}
\bigsetdef{U\in\cU}{\spec U \subset\kappa(K^c\cap\R)\cup\{1\}}. 
\end{split}
\end{equation}
Since $K$ is compact the set $\kappa(K^c\cap\R)\cup\{1\}$ is open.
Consequently 
\[
\bigsetdef{U\in\cU}{\spec U
\subset\kappa(K^c\cap\R)\cup\{1\}}
\]
is open and since $\Kappa$ is a
homeomorphism we reach the first conclusion.

The proof for ${\varrho}_\mathrm{ess}(T)$ instead of ${\varrho}(T)$
proceeds along the same lines.
\end{proof}

\begin{lemma}\label{S2.6} Let $K\subset\C$ be a compact set and let
$\Omega:=\bigsetdef{T\in\cC^{\sa}}{K\subset{\varrho}(T)}$ be equipped with
the gap topology. Then the map $R:K\times\Omega\to \cB, (\gl,T)\mapsto (T-
\gl)^{-1}$ is continuous. 
\end{lemma}

\begin{proof} For $(\gl,T)\in K\times \Omega$ we have
\begin{equation}\begin{split}
    R(\gl,T)=(T-\gl)^{-1}&=(I-(i+\gl)(T+i)^{-1})^{-1}(T+i)^{-1}\\
                       &=:F(\gl,(T+i)^{-1})=:F\circ G(\gl,T).
\end{split}
\end{equation}
In view of Theorem \plref{S1.1} the map
\begin{equation}
   \begin{split}
      G:&K\times \Omega\longrightarrow
K\times\bigsetdef{S\in\cB^{\sa}}{(K+i)^{-1}\subset{\varrho}(S)}\\        
&(\gl,T)\mapsto (\gl,(T+i)^{-1})
   \end{split}
\end{equation}
is continuous. Furthermore, the map
\begin{equation} \begin{matrix}
     F:&K\times\bigsetdef{S\in\cB}{(K+i)^{-1}\subset{\rho}(S)}
&\longrightarrow &\cB\\       
\ &(\gl,S)&\mapsto &(I-(i+\gl)S)^{-1}S
\end{matrix}
\end{equation}
is continuous. This proves the assertion.
\end{proof}

\begin{lemma}\label{S2.7} Let $a<b$ be real numbers. Then the
set 
\[
\Omega_{a,b}:=\bigsetdef{T\in\cC^{\sa}}{a,b\not\in\spec T}
\]
is open in
the gap topology and the map
\[\Omega_{a,b}\to \cB,\qquad T\mapsto 1_{[a,b]}(T)\]
is continuous.
\end{lemma}
\begin{proof} That $\Omega_{a,b}$ is open follows from Proposition
\plref{S2.1}. Next, denote by $\Gamma$ the circle of radius $(b-a)/2$ and
centre $(a+b)/2$. Then \begin{equation}
    1_{[a,b]}(T)=\frac{1}{2\pi i}\int_\Gamma (\gl-T)^{-1}d\gl.
\end{equation}
The assertion now follows from Lemma \plref{S2.6}.
\end{proof}

We collect what we have so far:

\begin{prop}\label{p:phillips-cont} 
Fix $T_0\in\cC\!\cF^{\sa}$. (a) Then there is a positive number $a$ and an open
neighbourhood $\cN \< \cC\!\cF^{\sa}$ of $T_0$ in the gap topology such that
the map 
\[
   \cN\to \cB,\qquad T\mapsto 1_{[-a,a]}(T)
\]
is continuous and finite-rank projection-valued, and hence 
$T\mapsto T1_{[-a,a]}(T)$ is also continuous. (We may as well assume
the rank to be constant).

\noi (b) If $-a\le c < d \le a$ are points such that $c,d\not\in \spec(T)$
for all $T\in \cN$ then the map $T\mapsto 1_{[c,d]}(T)$ is continuous on
$\cN$ and has finite rank on $\cN$. Of course, on any connected subset of
$\cN$ this rank is constant.
\end{prop}

\begin{proof} 
$T_0\in\cC\!\cF^{\sa}$ is equivalent to $0\not\in\specess(T_0)$. Thus either
$0\not\in\spec T_0$ or $0$ is an isolated point of $\spec T_0$ and an
eigenvalue of finite multiplicity. Hence there is an $a>0$ such that $\spec
T\cap [-a,a]\subset\{0\}$. By Lemma \ref{S2.5} the set
\begin{equation}
  \cN:=\bigsetdef{T\in\cC^{\sa}}{[-a,a]\subset {\varrho}_\mathrm{ess}(T), 
\tand \pm a\not\in\spec(T)} 
\end{equation}
is open in the gap topology and the map $T\mapsto 1_{[-a,a]}(T)$ is
continuous by Lemma \plref{S2.7}. 
Moreover, $\cN\<\cC\!\cF^{\sa}$ and $1_{[-a,a]}(T)$ is of finite rank.
This follows from the fact
that $[-a,a]\subset {\varrho}_\mathrm{ess}(T)$. This proves (a).
Now (b) follows from Lemma \plref{S2.7}.
\end{proof}

\begin{remark}\label{r:phillips-cont}
The preceding proposition is a precise copy of the corresponding result for
norm-continuous curves of bounded self-adjoint Fredholm operators. It
explains why, after all, spectral flow of gap-topology continuous
curves of (possibly unbounded) self-adjoint Fredholm operators can be
defined in precisely the same way as in the bounded case and with the same
properties. In substance, the proposition
was announced in \cite[p. 140]{BooWoj:EBP} without proof but with reference
to
\cite[IV.3.5]{Kat:PTL} (the continuity of a finite system of eigenvalues).
\end{remark}

Now we proceed exactly as in \cite[p. 462]{Phi:SAF}. We strive for almost 
literal repetition to emphasize the analogy (and the differences wherever
they occur) between the bounded and the unbounded case.

First a notation:
If $E$ is a finite-rank spectral projection for a self-adjoint operator
$T$, let $E^\ge$ denote the projection on the subspace of $E(H)$ spanned by
those eigenvectors for $T$ in $E(H)$ having non-negative eigenvalues.

\begin{dfn}\label{d:phillips-old}
Let $f:[0,1]\to \cC\!\cF^{\sa}(H)$ be a continuous path. By compactness and the
previous proposition, choose a partition, $\{0=t_0<t_1<\dots<t_n=1\}$ of
the interval and positive real numbers $\eps_j$, $j=1,\dots, n$ such that for
each $j=1,2,\dots,n$ the function $t\mapsto E_j(t):= 1_{[-\eps_j,\eps_j]}
\bigl(f(t)\bigr)$ is continuous and of finite rank on $[t_{j-1}, t_j]$. We
redefine the {\em spectral flow} of $f$, $\SF(f)$ to be
\[
\sum_{j=1}^n \Bigl( \dim\bigl(E_j^\ge(t_j)\bigr) 
- \dim\bigl(E_j^\ge(t_{j-1})\bigr) \Bigr).
\]
\end{dfn}

By definition, spectral flow is path additive when defined this way,
and we obtain in exactly the same way as in \cite{Phi:SAF}:

\begin{prop}\label{p:phillips-2}
Spectral flow is well defined; that is, it depends only on the continuous
mapping $f:[0,1]\to\cC\!\cF^{\sa}$\,.
\end{prop}

Propositions \plref{p:phillips-cont} and \plref{p:phillips-2} show that 
pathological examples like piecewise linear curves of self-adjoint 
unbounded Fredholm 
operators with infinitely fast oscillating spectrum and hence without
well-defined spectral flow are excluded; more precisely, they
cannot be continuous in 
the gap topology.

\begin{example}\label{x:fuglede-curve}
Let $H$ be a separable Hilbert space and $\{e_k\}_{k\in\N}$ be a complete
orthonormal system. Consider the multiplication operator 
which is defined by
\[
T_0:\dom(T_0)\to H, \quad \sum_k a_k e_k \mapsto \sum_k k a_k e_k 
\]
with $\dom(T_0)= \bigsetdef{\sum_k a_k e_k} {\sum_k k^2 |a_k|^2 <
+\infty}$. Then $T_0$ is self-adjoint and invertible and so $T_0\in
\cC\!\cF^{\sa}$. Set  \[
P_n:H\to H, \qquad e_k\mapsto
\begin{cases}
k e_k, & \text{ if   $k=n$,}\\
0, & \text{ otherwise}\ .
\end{cases}
\]
Then the sequence of unbounded self-adjoint Fredholm operators $\{T_{n}
:= 
T_0-2P_n\}_{n\in\N}$ converges to $T_0$ for $n\to\infty$ in the gap
topology. 
To see this, we apply Theorem \plref{S1.1}a and get
\begin{equation}\label{e:convergence}
\gamma(T_{n},T_0) = \Bigl\|(T_{n}+iI)\ii - (T_{0}+iI)\ii\Bigr\| =
\Bigl|\frac 1{i-n} - \frac 1{i+n}\Bigr| \to 0 \text{ for $n\to\infty$}. 
\end{equation}
For the Riesz transformation we note, however, that
\[
\|F_{T_{n}} e_n - F_{T_0} e_n\| = \Bigl|\frac {2n}{\sqrt{1+n^2}}\Bigr|
\to 2 \text{ for $n\to\infty$}.
\]
This is the aforementioned Fuglede example. Clearly the full spectrum (i.e.
the parts which are increasingly remote from 0) does
not change continuously for $n\to\infty$. The corresponding linear
interpolations $(1-t)T_{n}+tT_{n+1}$ all belong to $\cC\!\cF^{\sa}$ and
have rapidly oscillating spectrum also near $0$, hence the piecewise
linear curve can not be continuous in the gap topology by the previous
proposition; and it is not, as clearly seen by Theorem \plref{S1.1}a. We
find e.g.
\[\begin{split}
\gamma\bigl(\12 T_{n} + \12 T_{n+1}, T_0\bigr)
   &\ge \Bigl\|\Bigl(
\bigl(\12 T_{n} + \12 T_{n+1}+i\bigr)^{-1}-\bigl(T_0+i\bigr)^{-1}
\Bigr)e_n\Bigr\|\\
   &= \Bigl|\frac 1i - \frac 1{i+n}\Bigr|\to 1 \text{ for $n\to\infty$}.
\end{split}
\]

The example also shows that it is unlikely that
the Cayley image ${_{\cF}\!\cU}_{\!\inj}$ of
$\cC\!\cF^{\sa}$ can be retracted to the subspace where $1$ does not
belong to the spectrum at all (that is, the image of $\cF^{\sa}$ 
in ${_\cF\cU}_{\!\inj}$). Differently put, it shows 
that the eigenvalues of the
Cayley transforms flip around  $+1$ in the same way that the eigenvalues of
the operators in $\cC\!\cF^{\sa}$ flip around $\pm\infty$. More precisely,
consider the sequence of Cayley transforms
$U_{n}:=\Kappa(T_{n}) \in {_{\cF}\!\cU}_{\!\inj}$. The spectrum of $U_{n}$
consists of discrete eigenvalues which all lie in the lower half-plane
except one in the upper-half plane with a corresponding hole in the
lower half-plane sequence, plus the accumulation point 1 where $U_{n}-I$
is injective,  but not invertible. The same is true for
$U_{0}:=\Kappa(T_{0})$, but now having  all eigenvalues in the lower half
plane. By \eqref{e:convergence} the sequence $\{U_{n}\}_{n\in\N}$
converges to $U_0$ in ${_{\cF}\!\cU}_{\!\inj}$. We see that  the eigenvalues of
the sequence flip from the upper half-plane to the lower half-plane
close to +1 without actually crossing +1. It seems, however,
unlikely that there is a continuous path from $U_1$ to  $U_0$ which 
avoids any crossing.

Note that the linear path from $T_0$ to $T_1$ is continuous and 
has $\SF$ equal to -1. The corresponding curve from $U_0$ to $U_1$ has
one crossing at -1 from the lower half-plane to the upper one.

\end{example}

\medskip

\begin{remark}\label{r:homotopy}
So far we have established that spectral flow based on the
approach in \cite{Phi:SAF}, i.e., Definition \plref{d:phillips-old},
is well defined for gap continuous paths of self-adjoint Fredholm
operators. To do this we have
repeatedly used the local continuity proposition 
(Proposition \plref{p:phillips-cont}) for
continuous families in the gap topology. The surprising fact is that this
same local continuity proposition suffices to prove the homotopy
invariance. Initially, this may sound a little counter-intuitive since we
admit varying domains for our operators and therefore might not
expect nice parametrizations of the spectrum for these perturbations.

Of course it would suffice to show that Definition \plref{d:phillips-old}
coincides with the previous definition based on the
Cayley transform and the winding number (Definition \plref{S2.2}). 
Then, the homotopy invariance of Definition \plref{d:phillips-old}
would follow from Proposition \plref{S2.4} which is
based on general topological arguments. We prefer, however, to emphasize 
the existence of a self-contained proof based only on Definition \plref{S2.2} and 
Proposition \plref{p:phillips-cont}.
\end{remark}

\begin{prop}\label{p:phillips-3}
Spectral flow as defined in Definition \plref{d:phillips-old} is homotopy 
invariant.
\end{prop}

\begin{proof} As in \cite{Phi:SAF}.
\end{proof}

As a direct consequence of Proposition \plref{p:wind-def} we obtain:

\begin{prop}\label{p:coincide}
Spectral flows as defined in Definitions \plref{S2.2} and
\plref{d:phillips-old} coincide.
\end{prop}

\begin{remark}
In spite of the density of $\cF^{\sa}$ in $\cC\!\cF^{\sa}$ (Remark
\plref{r:open})
not any gap continuous path in $\cC\!\cF^{\sa}$ with endpoints in $\cF^{\sa}$ 
can be continuously deformed into an operator norm continuous path in
$\cF^{\sa}$.
One reason is that the one space is connected, but not the other by
Theorem \plref{t:connected}a.
\end{remark}

\section{Operator curves on manifolds with boundary}

In low-dimensional topology and quantum field theory, various examples of
operator curves appear which take their departure in a symmetric elliptic
differential operator of first order (usually an operator of Dirac type) on
a fixed compact Riemannian smooth manifold $M$ with boundary $\Sigma$.
Posing a suitable well-posed boundary value problem provides for a
nicely spaced discrete spectrum near 0. Then, varying the coefficients of
the differential operator and the imposed boundary condition suggests the
use of the powerful topological concept of spectral flow. In this Section
we show under which conditions the curves of the induced self-adjoint
$L^2$-extensions become continuous curves in $\cC\!\cF^{\sa}(L^2(M;E))$ in the
gap topology such that their spectral flow is well defined and truly
homotopy invariant.

\subsection{Notation and basic facts}
We fix the notation and recall basic facts, partially following
\cite{BooWoj:EBP} and \cite{Gru:TEP}. 

Let $D: \cinf{M;E}\to\cinf{M;E}$ be an elliptic symmetric (i.e., formally
self-adjoint) first order differential operator on $M$ acting on sections
of a Hermitian vector bundle $E$. Different from the case of closed
manifolds, now $D$ is no longer essentially self-adjoint and $\ker D$ is
infinite dimensional and varies with the regularity of the underlying
Sobolev space.  Among the
many extensions of $D$ to a closed operator in $L^2(M;E)$ we recall first
the definition of the {\em minimal} and the {\em maximal} closed extension
with
\begin{align*}
\dom(D^{{\rm min}}) &= \overline
{\bigsetdef{u\in \Ci(M;E)}{\supp u \< M\setminus \Sigma}}^{H^1(M;E)}
\tand\\
\dom(D^{{\rm max}}) &= \bigsetdef{u\in L^2(M;E)}{D u \in L^2(M;E)}.
\end{align*}

Now we make three basic (mutually related) assumptions:
\begin{asss}\label{a:dirac}
(1) The operator $D$ takes the form
\begin{equation}\label{e:product}
     D|_U=\sigma(y,\tau)\bigl(\frac{\dpa}{\dpa \tau}+A_\tau +
B_\tau\bigr)
\end{equation}
in a bi-collar $U={\Xi}\times[-\eps,\eps]$ of any hypersurface
$\Xi \< M\setminus \Sigma$, and a similar form in a collar of $\Sigma$,
where 
\begin{equation}
\sigma(\cdot,\tau),A_\tau, B_\tau: \Ci(\Xi_\tau;E|_{\Xi_\tau})
\too \Ci(\Xi_\tau;E|_{\Xi_\tau}) \label{final-3.2}
\end{equation} 
are a unitary bundle morphism; a symmetric elliptic differential operator
of first order; and a skew-symmetric bundle morphism, respectively, with 
\begin{equation}
\sigma(\cdot,\tau)^2=-I,\quad 
\sigma(\cdot,\tau) A_\tau = -A_\tau \sigma(\cdot,\tau), \tand
\sigma(\cdot,\tau) B_\tau = B_\tau \sigma(\cdot,\tau). \label{final-3.3}
\end{equation} 
Here $\tau$ denotes the normal variable and $\Xi_\tau$ a hypersurface
parallel to $\Xi$ in a distance $\tau$.

\smallskip

\noi (2) The operator $D$ satisfies the (weak) Unique Continuation Property
\begin{equation}\label{e:ucp}
\ker D^{{\rm max}} \cap \dom(D^{{\rm min}}) = \{0\}.
\end{equation}

\smallskip

\noi (3) The operator $D$ can be continued to an {\em invertible} 
elliptic differential operator $\wt D$ on a closed smooth Riemannian
manifold $\wt M$ which contains $M$ and acting on sections in a smooth
Hermitian bundle $\wt E$ which is a smooth continuation of $E$ over the
whole of $\wt M$; in particular, $\wt M$ is partitioned by 
${\Sigma}$ so that we have $\wt M=M_-\cup_{\Sigma} M_+$ with $M_+=M$,  $M_-
\cap M_+=\partial M_\pm = \Sigma$.

\end{asss}

\begin{remark}\label{r:ass}
All (compatible) Dirac operators satisfy Assumption (1) (see e.g.
\cite{Boo:UCP} or \cite{Gru:TEP}). Then Assumptions (2) and the sharper (3)
follow by \cite[Chapters 8, 9]{BooWoj:EBP}. 
\end{remark}

\bigskip

Let $\tilde \varrho, \varrho^\pm$ denote the trace maps from
$\cinf{\wt M;E}, \cinf{M_\pm;E}$ to $\cinf{{\Sigma};E|_{\Sigma}}$. (We
write $E$ also for $\wt E$ and $\wt E|_{M_-}$). Furthermore, $r^\pm$
denotes restriction to $M_\pm$ and $e^\pm$ denotes extension by $0$ from
$M_\pm$ to $\tilde M$. 

Under the fundamental Assumption (3) it is
well known that the \emph{Poisson operator} $K$ is given by
\begin{equation}\label{e:poisson}
     K:=r^+\tilde D^{-1}\tilde\varrho^*\sigma.
\end{equation}
The Poisson operator $K$ extends to a bounded mapping of $H^s(E|_{\Sigma})$
onto 
\[
Z^{s+1/2} =\bigsetdef{u\in H^{s+1/2}(M_+;E) }{Du=0 \text{ in the interior
of $M_+$}}
\]
and provides a left inverse for $\varrho^+:Z^{s+1/2}\to
H^s(E{|_{\Sigma}})$. Note that by the ellipticity of $D$, the trace map
$\varrho^+$ can be extended to $Z^{s+1/2}$ for all real $s$
(cf. \cite[Theorem 12.4]{BooWoj:EBP}). 

The \Calderon\ projector is then given by
\begin{equation}\label{e:calderon}
   P_+= \varrho^+ K.
\end{equation}
It is a pseudodifferential projection (idempotent). By definition, its
extension to $H^s(E_{|{\Sigma}})$ has the {\em Cauchy data space}
$\varrho^+ (Z^{s+1/2})$ as its range. Without loss of generality we can
assume that the extension of $P_+$ to $L^2(E|_{{\Sigma}})$ is orthogonal
(see \cite[Lemma 12.8]{BooWoj:EBP}).

\subsection{{\em Well-posed} boundary problems}
To obtain self-adjoint Fredholm extensions of $D$ in $L^2(M_+;E)$ we must
impose suitable boundary conditions. 

\begin{dfn}\label{d:grass}
The self-adjoint Fredholm Grassmannian of $D$ is defined by
\begin{multline*}
{\rm Gr}^{\sa}(D):= 
\bigl\{P\, \text{pseudodifferential projection}\mid 
P^*=P,\, P=\sigma_0(I-P)\sigma_0^*,\\
 \text{and $PP_+:\mathrm{range} P_+\to\mathrm{range} P$ Fredholm}\bigr\},
\end{multline*}
where $\sigma_0:E|_\Sigma\to E|_\Sigma$ denotes the unitary bundle morphism
over the boundary according to Assumption (1). The topology is given
by the operator norm.
\end{dfn}

It is well known (see e.g. \cite[Appendix B]{DouWoj:ALE}) 
that ${\rm Gr}^{\sa}(D)$ is connected
with the higher homotopy groups given by Bott periodicity.

\begin{remark}\label{r:aps}
At $\Sigma$, the ``tangential operator" $A_0$ defines a spectral projection
$\Pi_\ge$ of $L^2(\Sigma;E|_\Sigma)$ onto the subspace spanned by the
eigensections of $A_0$ for non-negative eigenvalues, the {\em
Atiyah-Patodi-Singer projection}. If $A_0$ is invertible, then
$\Pi_\ge=\Pi_>$  belongs to ${\rm Gr}^{\sa}(D)$. If $A_0$ is not invertible,
then one adds to $\Pi_>$ a projection onto a Lagrangian subspace of $\ker
A_0$ (relative to $\sigma_0$) to obtain an element in ${\rm Gr}^{\sa}(D)$. 
\end{remark}

We recall the main result of the analysis of well-posed boundary problems
(see e.g. \cite[Corollary 19.2, Theorem 19.5, and Proposition 20.3]{BooWoj:EBP}):

\begin{theorem}\label{t:grass} 
(a)
Each $P\in {\rm Gr}^{\sa}(D)$ defines a self-adjoint extension $D_P$ in
$L^2(M;E)$ with compact resolvent by
\[
\dom(D_P):= \bigsetdef {u\in H^1(M;E)}{P(u|_\Sigma)=0}.
\]

\noi (b) The \Calderon\ extension $D_{P_+}$ is invertible. In
fact, the inverse of $D_{P_+}$ can be expressed in terms of $\tilde D^{-1}$
and the Poisson operator: 
\begin{equation}
   D_{P_+}^{-1}=r^+\tilde D^{-1}e^+-KP_+\tilde\varrho \tilde D^{-1}e^+.
\end{equation}

\noi (c) The operator $D_P$ is invertible if and only if the {\em boundary
integral} 
\[
P\circ P_+:\mathrm{range} P_+\to \mathrm{range} P
\]
is invertible. Denote by $\tilde
Q_P$ its inverse and put $Q_P:=\tilde Q_PP$.
Then 
\begin{equation}\label{e:inverse}\begin{split}
    D_P^{-1}&=D_{P_+}^{-1}-KQ_P\varrho^+D_{P_+}^{-1}\\
       &=(I-KQ_P\varrho^+)\bigl(r^+ \tilde
D^{-1}e^+-KP_+\tilde\varrho\tilde
D^{-1}e^+\bigr). \end{split}
\end{equation}
\end{theorem}

\begin{lemma}\label{l:grass}
Let $H$ be a Hilbert space. For an invertible pair $(P,R)$ of 
orthogonal projections let $\wt Q(P,R)$ denote the inverse of
$PR:\mathrm{range} R \to \mathrm{range} P$ and put 
\[
Q(P,R):=\wt Q(P,R)P. 
\]
Then the map
\[
(P,R) \mapsto Q(P,R) \in \cB(H)
\]
is continuous in the operator norm.
\end{lemma}

\begin{proof} $(P,R)$ is an invertible pair if and only if
\[
T(P,R) := PR+(I-P)(I-R)
\]
is an invertible operator. Obviously, $(P,R)\mapsto T(P,R)$
is continuous on the set of invertible pairs. From
\[
T(P,R)R=PR=PT(P,R), \quad T(P,R)(I-R)=(I-P)T(P,R)
\]
we infer
\[
RT(P,R)\ii = T(P,R)\ii P, \quad (I-R)T(P,R)\ii = T(P,R)\ii(I-P)
\]
and so $Q(P,R)=T(P,R)\ii P$, and we reach the conclusion.
\end{proof}

\begin{cor}\label{c:grass}
For fixed $D$ the mapping 
\[
{\rm Gr}^{\sa}(D) \ni P \mapsto D_P \in \cC\!\cF^{\sa}(L^2(M;E))
\]
is continuous from the operator norm to the gap metric.
\end{cor}

\begin{proof} It follows immediately from \eqref{e:inverse}, 
Lemma \plref{l:grass}, and Theorem \plref{S1.1}a (see also 
Remark \plref{r:metrics}a) that
\[
\bigl\{P\in {\rm Gr}^{\sa}(D)\mid (P,P_+) \text{ invertible}\}\quad \ni\quad P 
\quad \mapsto \quad D_P \quad \in \quad \cC\!\cF^{\sa}(L^2(M;E))
\]
is continuous. Now consider $P_0\in {\rm Gr}^{\sa}(D)$ such that 
$D_{P_0}$ is not invertible. Since $D_{P_0}\in  \cC\!\cF^{\sa}(L^2(M;E))$,
the operator $D_{P_0}+\eps = (D+\eps)_{P_0}$ is invertible
for any real $\eps >0$ small enough. Obviously $D+\eps$ also satisfies
Assumptions \plref{a:dirac}, (1)-(3) and its invertible extension,
$\widetilde{D+\eps}$, may be choosen as $\tilde D+\eps$ which depends
continuously
on $\eps$. 
In view of \eqref{e:calderon}
the \Calderon\ projector $P_+(D+\eps)$ depends continuously
on $\eps$ (see also Theorem \plref{t:dependence} below). Thus for
$\eps$ small enough we have $P_0\in {\rm Gr}^{\sa}(D+\eps)$ with 
$\bigl(P_0,P_+(D+\eps)\bigr)$ invertible and the above argument
shows that
$
P\mapsto (D+\eps)_P=D_P+\eps
$
is continuous at $P_0$\,. Since $\eps=\eps\cdot I$ is bounded, also
$P\mapsto D_P$ is continuous at $P_0$\,.
\end{proof}

\bigskip

\subsection{The variation of the operator $D$}
We now assume 
that $D$ depends on an additional parameter $s$. More precisely, 
let $(D_s)_{s\in X}$\,, $X$ a metric space, be a family of differential 
operators satisfying the Assumption \plref{a:dirac} (1). We assume
moreover that 
\begin{equation}\label{e:3.9}
\text{{\em in each local chart, the coefficients of $D_s$ depend continuously on} $s$}.
\end{equation}
In a collar $U\approx [0,\eps)\times\Sigma$ the operator $D_s|_U$ 
takes the form
\begin{equation}\label{e:3.10}
     D_s|_U = \sigma_s(y,\tau)\bigl(\frac{\dpa}{\dpa \tau}+A_{s,\tau}
 + B_{s,\tau}\bigr)
\end{equation}
with $\gs,\, A,\, B$ depending continuously on $s$ and smoothly
on $\tau$. By the very definition of smoothness on a manifold with
boundary we find extensions of $\gs,\, A,\, B$ to 
\[
(s,\tau)\in X\times [-\gd,\eps)
\]
for some $\gd > 0$, such that \eqref{final-3.2} and \eqref{final-3.3} 
are preserved and such that the operator
\begin{equation}\label{e:3.11}
{D_s}':=\begin{cases} D_s, &\text{ on $M$,}\\
   \sigma_s\bigl(\frac{\dpa}{\dpa \tau}+A_{s,\tau}  + B_{s,\tau}\bigr),
& \text{ on $[-\gd,\eps)\times\gS$,}
\end{cases} 
\end{equation}
is a first order elliptic differential operator on the manifold 
$M_\gd:= \bigl([-\gd,0]\times \gS\bigr)\cup_\gS M$.

We fix $s_0\in X$. We choose $\gd$ so small that the operator
\begin{equation}\label{e:3.12}
\bigl[t \sigma_{s_0}(y,-\gd)\bigl(\frac{\dpa}{\dpa \tau}+A_{s_0,-\gd}
 + B_{s_0,-\gd}\bigr) + (1-t){D_s}'\bigr]\Big\rvert_{[-\gd,\gd]\times \gS}
\end{equation}
is elliptic for all $t\in [0,1]$.

Next we choose a cut-off function $\varphi\in\Ci(\R)$ with
\begin{equation}\label{e:3.13}
\varphi(x)=\begin{cases}
1,\ &\qquad x\le -\frac23\gd,\\
0, \ &\qquad x\ge -\frac13\gd\, .
\end{cases} 
\end{equation}  
Then we put 
\begin{equation}\label{e:3.14}
{D_s}'' := \varphi 
\sigma_{s_0}(y,-\gd)
\bigl(\frac{\dpa}{\dpa \tau}+A_{s_0,-\gd}  + B_{s_0,-\gd}\bigr) 
+ (1-\varphi){D_s}'\,.
\end{equation}
Clearly, ${D_{s_0}}''$ is an elliptic differential operator on $M_\gd$
which satisfies assumption \eqref{e:product}. Moreover, in the 
collar $U'':=[-\gd,-\frac 23\gd]\times \gS$ of $\partial M_\gd$ we have
\begin{equation}\label{e:3.15}
{D_s}'' = \sigma_{s_0}(y,-\gd)
\bigl(\frac{\dpa}{\dpa \tau}+A_{s_0,-\gd}  + B_{s_0,-\gd}\bigr) 
=: \sigma'' \bigl(\frac{\dpa}{\dpa \tau}+A''  + B''\bigr) ,
\end{equation}
where $\gs'',\, A'',\, B''$ are {\em independent} of $s$ and $\tau$.

By construction, ${D_s}''$ preserves \eqref{e:3.9}. Hence there is an open
neighbourhood $X_0$ of $s_0$ such that for $s\in X_0$ the operator 
${D_s}''$ is elliptic.

For $\{{D_s}''\}_{s\in X_0}$ we now apply the construction of the 
invertible double of \cite[Chapter 9]{BooWoj:EBP}. In view of 
\eqref{e:3.15}, the invertible double will be a first order elliptic
differential operator on a closed manifold which depends continuously
on the parameter $s$.

Summing up we have proved

\begin{theorem}\label{t:3.8}
Let $M$ be a compact Riemannian manifold with boundary. Let 
$\{D_s\}_{s\in X}$ \,, $X$ a metric space, be a family of differential
operators satisfying Assumption \plref{a:dirac}.(1) and which
depends continuously on $s$ in the sense of \eqref{e:3.9}. Then for
each $s_0\in X$ there exist an open neighbourhood $X_0$ of $s_0$
and a continuous family $\{\wt{D_s}\}_{s\in X_0}$ 
of invertible elliptic differential operators 
$\wt{D_s}:\Ci(\wt M;\wt E)\to\Ci(\wt M;\wt E)$ with $\wt {D_s}|_M=D_s$\,.
Here $\wt M$ is a closed 
Riemannian manifold with $\wt M\> M$; and $\wt E\to\wt M$ a smooth Hermitian
vector bundle with $\wt E|_M=E$.
\end{theorem}

The continuity of $s\mapsto \wt{D_s}$ is understood in the sense of 
\eqref{e:3.9}. However, since $\wt M$ is closed this implies that 
$\{\wt{D_s}\}_{s\in X_0}$ is a graph continuous family of invertible
self-adjoint operators.

\begin{theorem}\label{t:dependence}
Under the assumptions of Theorem \plref{t:3.8} we have

\noi (a) The Poisson operator $K_s$ of $D_s$ depends continuously
on $s$.

\noi(b) The \Calderon\ projector $P_+(s)$ of $D_s$ depends
continuously on $s$.

\noi(c) The family 
\[
X\ni s\qquad \mapsto \qquad (D_s)_{P_+(s)} \in \cC\!\cF^{\sa}(L^2(M;E))
\]
is continuous.

\noi (d) Let $\{P_t\}_{t\in Y}$ be a norm-continuous path of orthogonal
projections in $L^2(\gS;E|_\gS)$\,. If
\[
P_t\in \bigcap_{s\in X}{\rm Gr}^{\sa}(D_s), \qquad t\in Y\,,
\]
then
\[
X\times Y \quad\ni \quad (s,t)\quad \mapsto \quad (D_s)_{P_t} \quad
\in\quad \cC\!\cF^{\sa}(L^2(M;E))
\]
is continuous.
\end{theorem}

\begin{proof}\footnote{\textsc{Added in Proof:} The proof of Theorem
    \plref{t:dependence} is incomplete. The continuous dependence of
 formula
\eqref{t:3.8} on all input data is obvious only if $K$ is viewed as
a map from $H^s(E|_{\Sigma})$ to $Z^{s+1/2}$ for $s<0$. In the critical
case $s=0$ additional considerations are necessary (for the Trace theorem
in the critical Sobolev case $s=0$ see e.g. \cite[Theorem 12.4]{BooWoj:EBP}).
Nevertheless, Theorem \plref{t:dependence} is correct, though 
the proof is more involved. For a perturbative approach
see Section 3 of B. Himpel, P. Kirk, and M. Lesch: \textit{
Calderon projector for the Hessian of the perturbed Chern-Simons function 
on a 3-manifold with boundary. } 
To appear in Proc. London Math. Soc.; math.GT/0302234}
(a) follows from Theorem \plref{t:3.8} and \eqref{e:poisson};
(b) follows from  Theorem \plref{t:3.8} and \eqref{e:calderon};
(c) follows from  Theorem \plref{t:3.8} and \eqref{e:poisson}.

\noi (d) Similarly as in the proof of Lemma \plref{l:grass} it suffices to prove
the claim for $(D_s)_{P_t}$ invertible. Now the assertion follows from
Lemma \plref{l:grass} and \eqref{e:inverse}.
\end{proof}

\begin{remark}\label{r:nicolaescu}
(a) By different methods, somewhat related results have been obtained
in \cite{BooFur:MIF} under the additional assumption of a fixed 
principal symbol of the family $\{D_s\}$ and a fixed boundary 
condition.

\noi (b) Corollary \plref{c:grass} for fixed $D$ and the preceding
Theorem \plref{t:dependence} for variation of $D$ yield a 
well-defined and homotopy invariant spectral flow by Proposition
\plref{S2.3}, resp. Propositions \plref{p:phillips-3}, \plref{p:coincide}.
The surprising facts are that 
\begin{enumerate}
\item gap continuity suffices to establish spectral flow and
\item gap continuity is obtainable from continuous variation of the
operator and the boundary condition without any restrictions and without
any need to fix the domains of the unbounded $L^2$-extensions by unitary
transformations.
\end{enumerate}
Roughly speaking, these constitute the differences between the present approach
and Nicolaescu's approach in \cite{Nic:SFO} which requires the continuity
of the Riesz map and to achieve that additional properties of the families
of boundary problems.

\noi (c) In some important applications in topology, families of
Dirac operators are considered on non-compact manifolds.
The $L^2$-extensions of these operators are self-adjoint Fredholm operators but
do not have a compact resolvent and therefore require a
light modification of our preceding arguments to establish the continuity
in the gap metric.
\end{remark}


\end{document}